\documentclass{article}
\usepackage{amsthm}

\begin{document}

Tsemo Aristide

PKFokam Institute of Excellence,

Yaounde, Cameroon.

tsemoaristide@hotmail.com

\bigskip
\bigskip

\centerline{\bf Globally   linearizable actions on topological manifolds.}

\bigskip
\bigskip

\centerline{\bf Abstract.}

\bigskip

{\it Let $M$ be a finite dimensional topological aspherical manifold whose universal cover is ${\bf R}^n$. In this paper, we study $Aff(M)$, the subgroup of the group of homeomorphisms of $M$, whose elements can be lifted to affine transformations of ${\bf R}^n$. We show that if $M$ is closed, the connected component $Aff(M)_0$ of $Aff(M)$ acts locally freely on $M$. We deduce that $Aff(M)_0$ is a solvable Lie group, and is nilpotent if $M$ is a polynomial manifold. We study the foliation defined by the orbits of $Aff(M)_0$ if $dim(Aff(M)_0)=dim(M)-1$.}

\bigskip
\bigskip

\centerline{\bf I. Introduction.}

\bigskip

Let $M$ be a topological $n$-dimensional manifold, and $G$ a topological group which acts on $M$. The action of $G$ on $M$ is linearizable at $x$, if there exists a neighbourhood $U$ of $x$ such that $U$ is stable by $G$, and the restriction of $G$ to $U$ is conjugated to a linear action. If $M$ is the Euclidean space ${\bf R}^n$, we say that the action is globally linearizable if we can choose $U={\bf R}^n$. Remark that an infinitesimal version of this problem has been studied by Poincare in his thesis. He has given conditions which ensure that an ordinary nonlinear  differential equation is equivalent to a linear differential equation. A generalization of this Poincare's problem is given by a finite set $L$ of differential equations which form a Lie algebra ${\cal L}$. Suppose that a group $G$ acts freely, properly and discontinuously by diffeomorphisms on ${\bf R}^n$, and the equations of $L$ are invariant by $G$. This implies that the quotient of ${\bf R}^n$ by $G$ is a manifold $M$, and there exists a Lie algebra ${\cal L}_G$ of $M$ whose lifts by the universal covering map of $M$ is ${\cal L}$. If  ${\cal L}$ is conjugated by a diffeomorphism of ${\bf R}^n$ to an affine system, or equivalently to a Lie algebra contained in the Lie algebra of the group of affine transformations of ${\bf R}^n$, we say that ${\cal L}$ is globally linearizable. Remark that, if $M$ is compact, ${\cal L}_G$ is the image of an infinitesimal action of a simply connected Lie group on $M$, whose action can be lifted to ${\bf R}^n$.

In this paper, we study a topological version of global linearizable actions. Let $M$ be an $n$-dimensional topological manifold whose universal cover is homeomorphic to ${\bf R}^n$. We denote by $\pi_1(M)$ the fundamental group of $M$. We consider the group $Aff(M)$ of homeomorphisms of $M$ which can be lifted to affine maps of ${\bf R}^n$. We characterize the subgroup $N(\pi_1(M))$ of homeomorphisms of ${\bf R}^n$ which are lifts of elements of $Aff(M)$ by showing that it is generated by its normal subgroup $\pi_1(M)$ and a subgroup of the group of affine transformations of ${\bf R}^n$. We adapt a result of Fried, Goldman and Hirsch  by showing that if $M$ is closed, the action of $\pi_1(M)$ on ${\bf R}^n$ is irreducible, that is, does not preserve any proper non empty affine subset. This enables to show that the action of the connected component $Aff(M)_0$ of $Aff(M)$ on $M$ is locally free. We deduce that if $M$ is a compact  $Aff(M)_0$ is solvable, and if $M$ is a compact polynomial manifold, $Aff(M)_0$ is nilpotent. We remark that if $Aff(M)_0$ acts transitively on $M$, then $M$ is a solvmanifold endowed with an affine structure. Finally, we  determine the topology of $M$, if $M$ is differentiable and $dim(Aff(M)_0)=dim(M)-1$.
\bigskip
\bigskip

{\bf 2. General properties of the group of affine transformations of a topological manifold.}

\bigskip

In this section, we are going to define affine transformations of topological manifolds and study some of their general properties.

\bigskip

{\bf Definition 2.1.}
A group of homeomorphisms $G$ of the Euclidean space ${\bf R}^n$ acts freely and properly discontinuously on ${\bf R}^n$ if and only if:

For every elements $x$ of ${\bf R}^n$, and $\gamma$ of $G$, $\gamma(x)=x$ implies that $\gamma=Id_{{\bf R}^n}$.

For every element $x$ of ${\bf R}^n$, there exists an open subset $U$ containing $x$ such that for every $\gamma$ in $G$, $\gamma(U)\cap U$ is not empty implies that $\gamma=Id_{{\bf R}^n}$.

\medskip

For such an action, the quotient of ${\bf R}^n$ by $G$ is a topological manifold $M$. We  denote by $p_M:{\bf R}^n\rightarrow M$ the quotient map. It is a covering map. We can define:

\medskip

{\bf Definition 2.2.}
An ${\bf R}^n$-aspherical manifold $M$, is a topological manifold $M$, such that there exists a group $G$  which acts continuously, freely, and properly discontinuously on ${\bf R}^n$, and the quotient of ${\bf R}^n$ by $G$ is $M$.

\medskip

 We can identify the fundamental group $\pi_1(M)$, of $M$, with $G$.

\bigskip

{\bf Definitions 2.3.}
Let $G$ be a group which acts on ${\bf R}^n$ by homeomorphisms.
An affine subset $V$ of ${\bf R}^n$ is stable by $G$ if for every $x$ in $V$, and every $g$ in $G$, $g(x)$ is an element of $V$.

The action of $G$ on ${\bf R}^n$ is irreducible if the only affine subsets of ${\bf R}^n$ stable by $G$ are the empty subset and ${\bf R}^n$. 

\medskip

The following fundamental proposition for this paper is adapted from   [3], theorem 2.2 p. 496.

\bigskip

{\bf Proposition 2.1.}
{\it Let $X$ be an $n$-dimensional contractible topological manifold, suppose that a group $G$ acts freely and properly discontinuously on $X$ and the quotient of $X$ by $G$ is compact. Let $Y$ be a non empty topological contractible submanifold of $X$. Suppose that for every element $g\in G, g(Y)\subset Y$, then $dim(Y)=n$.}

\medskip

{\bf Proof.}
The quotient $M$ of $X$ by $G$ (resp. $N$ of $Y$ by $G$)  are Eilenberg-McLane $K(G,1)$-spaces.  Up to a $2$-cover, we can assume that $M$ and $N$ are oriented, this implies that the singular homology group $H_n(M,{\bf Z})$ is not zero since $M$ is compact (see [1] p. 346). we deduce that $H_n(N,{\bf Z})$ is not zero. This implies that $dim(Y)= n$.

\bigskip

{\bf Corollary 2.1.}
{\it  Let $M$ be a compact ${\bf R}^n$-aspherical manifold, the action of $\pi_1(M)$ on ${\bf R}^n$ is irreducible.}

\medskip

{\bf Proof.}
Let $V$ be a non empty affine subspace such that for every $\gamma\in \pi_1(M)$, $\gamma(V)=V$. The proposition 2.1 implies that the dimension of the quotient of $V$ by $G$ is $n$. We deduce that $V={\bf R}^n$.

\bigskip

Let $M$ be an ${\bf R}^n$-aspherical  manifold, we can lift every homeomorphism $f$ of $M$ to an homeomorphism $\hat f$ of ${\bf R}^n$ such that the following diagram is commutative (see [1] Theorem 4.1):

$$
\matrix{ {\bf R}^n &{\buildrel{\hat f}\over{\longrightarrow}}& {\bf R}^n\cr p_M\downarrow &&\downarrow p_M\cr M&{\buildrel{f}\over{\longrightarrow}}& M}
$$

\medskip

{\bf Definition 2.4.}
An homeomorphism $f$ of the ${\bf R}^n$-aspherical manifold $M$ is an affine transformation, if there exists a lift $\hat f$ of $f$ which is an element of $Aff({\bf R}^n)$, the group of affine transformations of   ${\bf R}^n$. We denote by $Aff(M)$ the set of affine transformations of $M$.

\medskip

The following proposition can be proved straightforward by following the definitions:

\medskip

{\bf Proposition 2.2.}
{\it The set of affine transformations $Aff(M)$ is a group.}

\bigskip

Let $M$ be an ${\bf R}^n$-aspherical manifold, consider $N(\pi_1(M))$, the set of elements of the group of homeomorphisms of ${\bf R}^n$  such that, for every element $g\in N((\pi_1(M))$, there exists an element $f$ of $Aff(M)$ whose lift is $g$. Remark that $g$ is not necessarily an affine transformation of ${\bf R}^n$. For example, every element of $\pi_1(M)$ is an element of $N(\pi_1(M))$ since it is the lift of $Id_M$, the identity of $M$, and $Id_{{\bf R}^n}$ is the lift of $Id_M$.  We have the following proposition:

\medskip

{\bf Proposition 2.3.}
{\it Let $g$ be an element of $N(\pi_1(M))$, there exists an element $\hat f$ of $Aff({\bf R}^n)$, and an element $\gamma$ of $\pi_1(M)$ such that $g=\gamma\circ \hat f$.}

\medskip

{\bf Proof.}
There exists an element $f$ of $Aff(M)$ such that $g$ is the lift of $f$. Since $f$ is an element of $Aff(M)$, there exists $\hat f$ in $Aff({\bf R}^n)$ such that $\hat f$ is the lift of  $f$. Let $x$ be an element of ${\bf R}^n$, since $p_M(g(x))=p_M(\hat f(x))$, there exists $\gamma$ in $\pi_1(M)$ such that $\gamma(\hat f(x))=g(x)$. We are going to show that $g=\gamma\circ \hat f$.

The subset $V$ of elements of ${\bf R}^n$ on which $g$ and $\gamma\circ \hat f$ coincide is closed since $g$ and $\gamma\circ \hat f$ are continuous. Let $x$ be an element of $V$, there exists an open subset $U$ of ${\bf R}^n$ containing $g(x)$ such that, for every element $\alpha$ of $\pi_1(M)$, $\alpha(U)\cap U$ is not empty implies that $\alpha=Id_{{\bf R}^n}$. The open subset $g^{-1}(U)\cap (\gamma\circ \hat f)^{-1}(U)$ is not empty, since it contains $x$. For every element $y$ of $W=g^{-1}(U)\cap (\gamma\circ \hat f)^{-1}(U)$, $p_M(g(y))$ and $p_M(\hat f(x)))$ coincide. This implies that there exists $\alpha\in G$ such that $g(y)=\alpha( \hat f(y))$. We deduce that $\alpha (\hat f(y))$ and $\gamma(\hat f(y))$ are elements of $U$. We deduce that $\alpha\circ\gamma^{-1}(U)\cap U$ is not empty, since $(\alpha\circ\gamma^{-1})(\gamma(\hat f(y)))=\alpha(\hat f(y))$. This implies that $\alpha\circ\gamma^{-1}=Id_{{\bf R}^n}$, henceforth $\alpha=\gamma$. We deduce that $\gamma\circ \hat f$ and $g$ coincide on $W$. This implies that $V$ is open and closed. We deduce that $V$ is ${\bf R}^n$ since ${\bf R}^n$ is connected.

\bigskip

Let $G$ be a group of homeomorphisms of the topological space $X$. For every compact subset $K$ of $X$, and any open subset $U$ of $X$, we denote by $C_{K,U}$ the subset of $G$ such that for every $g\in C_{K,U}$, $g(K)$ is contained in $U$. The family $C_{K,U}$ is a basis of a topology defined on $G$ called the compact-open topology. We suppose that $Aff(M)$ and $N(\pi_1(M))$ are endowed with the compact-open topology. Remark that the compact-topology endows $\pi_1(M)$ with the structure of a discrete set:

\medskip

{\bf Proposition 2.4.}
{\it The compact-open topology endows $\pi_1(M)$ with the structure of a discrete set.}

\medskip

{\bf Proof.}
Let $x$ be an element of ${\bf R}^n$, there exists an open subset $U$ of ${\bf R}^n$ containing $x$ such that for every $\gamma\in \pi_1(M)$, $\gamma(U)\cap U$ is not empty if and only if $\gamma$ is the identity of ${\bf R}^n$. Let $\gamma$ be an element of $\pi_1(M)$. Let $C_{x,\gamma(U)}$ the subset of elements of $\pi_1(M)$ such that for every $\alpha$ of $C_{x,\gamma(U)}$, $\alpha(x)\subset \gamma(U)$. The subset $C_{x,\gamma(U)}$ is a neighbourhood of $\gamma$ for the compact-open topology defined on $\pi_1(M)$. If $\alpha$ is an element of $C_{x,\gamma(U)}$, the fact that $\alpha(x)\in \gamma(U)$, implies that $\gamma^{-1}(\alpha(x))$ is an element of $U$. We deduce that $(\gamma^{-1}\circ \alpha)(U)\cap U$ is not empty. The definition of $U$ implies  that $\gamma^{-1}\circ\alpha=Id_{{\bf R}^n}$, and $\gamma=\alpha$. We deduce that $C_{x,\gamma(U)}$ contains only $\gamma$, and $\pi_1(M)$ is a discrete subgroup for the compact-open topology.

\bigskip

{\bf Proposition 2.5.}
{\it An element of $N(\pi_1(M))$ which lifts the identity is contained in $\pi_1(M)$.}

\medskip

{\bf Proof.}
Let $f$ be an element of $N(\pi_1(M))$ above the identity, and $x$ an element of ${\bf R}^n$. There exists an element $\gamma$ of $\pi_1(M)$ such that $\gamma(x)=f(x)$, since $p_M(x)=p_M(f(x))$. Let $V$ be the subset of ${\bf R}^n$, such that for every $y\in V$, $f(x)=\gamma(x)$. The subset $V$ is closed since $\gamma$ and $f$ are continuous.  Since the action of $\pi_1(M)$ is  free and properly discontinuous, there exists  an open subset $U$ of ${\bf R}^n$ containing $f(x)$ such that for every element $\alpha$ in $\pi_1(M)$, $\alpha(U)\cap U$ is not empty implies that $\alpha$ is the identity. The open subset $W=\gamma^{-1}(U)\cap f^{-1}(U)$ is not empty since it contains $x$. Let $y\in W$, there exists $\alpha\in\pi_1(M)$ such that $f(y)=\alpha(y)$. This implies that, $\alpha^{-1}(f(y))=y$, and $\gamma(\alpha^{-1}(f(y)))=\gamma(y)$ is an element of $U$. We deduce that $(\gamma\circ\alpha^{-1})(U)\cap U$ is not empty. This implies that $\alpha=\gamma$, and henceforth $f=\gamma$.

\medskip

{\bf Proposition 2.6.}
{\it The group $\pi_1(M)$ is a normal subgroup of $N(\pi_1(M))$, and $Aff(M)$ is isomorphic to the quotient of $N(\pi_1(M))$ by $\pi_1(M)$.}

\medskip

{\bf Proof.}
Let $\gamma$ be an element of $\pi_1(M)$, and $f$ an element of $N(\pi_1(M))$, $f\circ \gamma\circ f^{-1}$ lifts the identity. We deduce that it is an element of $\pi_1(M)$, and $\pi_1(M)$ is normal in $N(\pi_1(M))$.

Let $P_M:N(\pi_1(M))\rightarrow Aff(M)$ the map which assigns to the element $\hat f$ of $N(\pi_1(M)))$, and element $f$ of $Aff(M)$ such that $\hat f$ lifts $f$. The map $f$ is well-defined. Suppose that $\hat f$ is the lift of the elements $f$ and $g$ of $Aff(M)$. Let $x$ an element of $M$, and $\hat x$ an element of $p^{-1}_M(x)$, we have $f(x)=g(x)=p_M(\hat f(\hat x))$.

The map $P_M$ is surjective, since every element of $Aff(M)$ can be lifted to $\hat M$.

The map $P_M$ is a morphism of groups. The proposition 2.5 implies that $P_M(\hat f)$ is the identity if and only if $\hat f$ is contained in $\pi_1(M)$. This is equivalent to saying that the kernel of $P_M$ is $\pi_1(M)$.

\bigskip

{\bf 3. The commutator of $\pi_1(M)$ and the action of $Aff(M)_0$.}

\bigskip

Let $Aff(M)_0$ be the connected component of $Aff(M)$, and $N(\pi_1(M))_0$ the connected component of $N(\pi_1(M))$. We denote by $Com(\pi_1(M))$, the subgroup of the group of homeomorphisms of ${\bf R}^n$ which commute with $\pi_1(M)$.  We endow $N(\pi_1(M))_0$ and $Aff(M)_0$ with the compact-open topology.

 We have the following result:

\medskip

{\bf Proposition 3.1.}
{\it The group $N(\pi_1(M))_0$ is contained in $Com(\pi_1(M))$. The image of the restriction of $P_M$ to $N(\pi_1(M))_0$ is $Aff(M)_0$.}

\medskip

{\bf Proof.}
Let $g$ an element of $N(\pi_1(M))_0$, and  $I=[0,1]$ the closed interval of ${\bf R}$, and  $f:I\rightarrow N(\pi_1(M))_0$ a continuous map such that $f(0)=Id_{{\bf R}^n}$ and $f(1)=g$. For every element $\gamma$ of $\pi_1(M)$, the map defined by $g_\gamma(t)=f(t)\circ\gamma\circ f(t)^{-1}$ is a continuous map contained in $\pi_1(M)$ since $\pi_1(M)$ is a normal subgroup of $N(\pi_1(M))_0$. Since $\pi_1(M)$ is discrete for the compact open topology, we deduce that $g_\gamma$ is constant and for every $t\in [0,1]$, $g_\gamma(t)=g_\gamma(0)=\gamma$. We deduce that every element of $N(\pi_1(M))_0$ commutes with $\pi_1(M)$.

The quotient map $P_M:N((\pi_1(M))\rightarrow Aff(M)$ is a covering map since its kernel is discrete. We deduce that we can lift any path through the identity  of $M$ contained in $Aff(M)_0$ to a path through the identity of ${\bf R}^n$ contained in $N(\pi_1(M))$. This implies that $P_M(N(\pi_1(M))_0)=Aff(M)_0$.

\bigskip

Let $Aff({\bf R}^n,M)$ be the intersection $N(\pi_1(M))\cap Aff({\bf R}^n)$. The restriction  of $P^A_M$ of $P_M$ to $Aff({\bf R}^n,M)$ is surjective, and its kernel is $Aff({\bf R}^n,M)\cap\pi_1(M)$. This implies that $P^A_M$ is a covering map since $\pi_1(M)\cap Aff({\bf R}^n,M)$ is discrete. Let $Aff({\bf R}^n,M)_0$ be the connected component of $Aff({\bf R}^n,M)$. Remark that since $Aff({\bf R}^n,M)_0$ is contained in $N(\pi_1(M))_0$, its elements commute with the elements of $\pi_1(M)$. We have the following result:

\bigskip

{\bf Proposition 3.2.}
{\it Suppose that $M$ is compact, then the group $Aff({\bf R}^n,M)_0$ acts freely on ${\bf R}^n$.}

\medskip

{\bf Proof.}
Let $f$ be an element of $Aff({\bf R}^n,M)_0$. Consider $Fix_{Aff}$ the set of affine subspaces of ${\bf R}^n$, such that for every element $V$ of $Fix_{Aff}$ and every $x$ in $V$, $f(x)=x$. Let $V$ be an element of $Fix_{Aff}$ which has a maximal dimension. Suppose that $dim(V)=p$ with $p<n$ and $V$ is not empty. The corollary 2.1 implies that there exist an element $x$ of $V$, and an element $\gamma$ of $\pi_1(M)$ such that $\gamma(x)$ is not contained in $V$.
Let $U$ be the direction of $V$, that is the vector subspace of dimension $p$ of ${\bf R}^n$ such that for every $y\in V$, there exists an element $u$ of $U$ such that $y=x+u$. We can write $f=(L,l)$ where $L$ is a linear map of ${\bf R}^n$, $l$ an element of ${\bf R}^n$, and for every $z$ in ${\bf R}^n$, $f(z)=L(z)+l$.

Write $\gamma(x)=x+v$. We have $f(\gamma(x))=\gamma(f(x))=\gamma(x)$. This implies that $f(\gamma(x))=L(x+v)+l=f(x)+L(v)=x+v$. Since $f(x)=x$, we deduce that $L(v)=v$. All the points of the affine subspace which contains $x$, and whose direction is the direct sum $U\oplus {\bf R}v$ are fixed by $f$. This implies that de dimension of $V$ is not maximal in $Fix_{Aff}$, contradiction. We deduce that either $V$ is ${\bf R}^n$ or $V$ is empty. This is equivalent to saying that $Aff({\bf R}^n,M)_0$ acts freely on ${\bf R}^n$.

\bigskip

Remark that if $f$ is an element of $Aff({\bf R}^n)$ which commutes with $\pi_1(M)$, it is the lift of an element of $Aff(M)$. The subgroup $Com({\bf R}^n,M)$ of $Aff({\bf R}^n)$ whose elements commute with the elements of $\pi_1(M)$ is a Lie group since it is a closed subgroup of $Aff({\bf R}^n)$. Since it contains $Aff({\bf R}^n,M)_0$, its connected component is the connected of $Aff({\bf R}^n,M)$. We deduce that $Aff({\bf R}^n,M)_0$ is a Lie group. We denote by $aff({\bf R}^n,M)$ the Lie algebra of $Aff({\bf R}^n,M)$.

\bigskip

{\bf Corollary 3.1.}
{\it Suppose that $M$ is compact, then $Aff({\bf R}^n,M)_0$ and $Aff(M)_0$ are  solvable Lie groups.}

\medskip

{\bf Proof.}
Let $KAN$ be the Iwasawa decomposition of $Aff({\bf R}^n,M)_0$. The group $K$ is a compact Lie group which acts freely on ${\bf R}^n$, we deduce that it is trivial. This implies that $Aff({\bf R}^n,M)_0$ is a solvable group. The restriction of $P_M$ to $Aff({\bf R}^n,M)_0$ is a covering map whose image is $Aff(M)_0$. This implies that $Aff(M)_0$ is solvable.

\bigskip

{\bf Corollary 3.2.}
{\it Suppose that $M$ is compact, then group $Aff(M)_0$ is a Lie group whose action on $M$ is locally free.}

\medskip

{\bf Proof.}
Let $x$ be an element of $M$, $\hat x$ an element of $p_M^{-1}(x)$. Consider an open subset $V$ of the identity of $Aff({\bf R}^n,M)_0$ such that the restriction of $P_M$ to $V$ is an homeomorphism onto its image. We can shrink $V$ and find a neighbourhood $U$ of $\hat x$ such that: the restriction of $p_M$ to $U$ is an homeomorphism onto its image, and for every element $\hat f$ of $V$, $\hat f(\hat x)$ is contained in $U$. Let $\hat f$ be an element of $V$, write $f=P_M(\hat f)$. Suppose that $f(x)=x$, this is equivalent to saying that  $p_M(\hat f(\hat x))=x$. Since $\hat f(\hat x)$ is contained in $U$, we deduce that $\hat f(x)=\hat x$, and $\hat f=Id_{{\bf R}^n}$ since $Aff({\bf R}^n,M)_0$ acts freely on ${\bf R}^n$. This implies that $f$ is the identity of $M$, and the action of $Aff(M)_0$ on $M$ is locally free.

\bigskip

We suppose now that the action of $\pi_1(M)$ on ${\bf R}^n$ is a differentiable action. This implies that $M$ is a differentiable manifold, we have the following result:

\bigskip

{\bf Corollary 3.3.}
{\it The orbits of $Aff(M)_0$ are the leaves of a foliation.}

\bigskip

We can also show:

\bigskip

{\bf Corollary 3.4.}
{\it Suppose that $Aff({\bf R}^n,M)_0=Com({\bf R}^n,M)_0$ is an algebraic group,  $Aff(M)_0$ is a nilpotent group.}

\bigskip

{\bf Proof.}
Let $U$ be the unipotent radical of $Aff({\bf R}^n,M)_0$, (see [6), there exists a reductive subgroup $R$ such that $Aff({\bf R}^n,M)_0$ is the semi-direct product of $U$ and $R$. Since any representation of a reductive subgroup has a fixed point, we deduce that $R$ is trivial and $Aff({\bf R}^n,M)_0$ is unipotent, and henceforth $Aff(M)_0$ is nilpotent.

\bigskip

A compact and complete polynomial manifold of dimension $n$ is a compact differentiable manifold which is  the quotient of ${\bf R}^n$ by a subgroup of polynomial automorphisms which act freely and properly discontinuously on ${\bf R}^n$. Since $M$ is compact, its fundamental group is generated by a finite set $S$ of polynomial automorphisms. This implies that $Com({\bf R}^n,M)$ is defined by a finite set of algebraic equations, and is an algebraic group. We have:

\bigskip

{\bf Corollary 3.5.}
{\it Suppose that $M$ is a compact and complete polynomial manifold, then $Aff(M)_0$ is nilpotent.}

\bigskip

{\bf Corollary 3.6.}
{\it Suppose that $M$ is a compact differentiable manifold, and $Aff(M)_0$ is a compact Lie group of dimension $l$. Then  $M$ is the total space of a bundle whose typical fibre is the torus ${\bf T}^l$ over an orbifold.}

\medskip

{\bf Proof.}
The fact that $Aff(M)_0$ is compact implies that it is commutative, since it is a solvable Lie group. We deduce that $M$ is endowed with the locally free action of ${\bf T}^l$. This  implies that that $M$ is the total space of a fibre bundle whose typical fibre is ${\bf T}^l$, and whose base space is an orbifold.

\bigskip

{\bf Proposition 3.3.}
{\it Suppose that $M$ is a compact differentiable manifold, which is the quotient of ${\bf R}^n$ by a free and properly discontinuous  action of its fundamental group. Moreover, suppose that the subset of elements of $Aff({\bf R}^n)$ which commute with $\pi_1(M)$ is not discrete, then the Euler number of $M$ is zero.}

\medskip

{\bf Proof.}
Under the hypothesis of the proposition, the Lie algebra $aff(M)$ of $Aff(M)$ is not zero, since it is isomorphic to the Lie algebra of $Com({\bf R}^n,M)$ which is not zero. Since $Aff(M)_0$ acts locally freely on $M$, any non zero element of $aff(M)$ does not vanish on $M$. We deduce that the Euler number of $M$ is zero.

\bigskip

{\bf Definitions 3.1.}
An affine manifold $(M,\nabla_M)$ is a manifold endowed with a connection whose curvature and torsion forms vanish identically. This is equivalent to endow $M$ with an atlas whose coordinates change are affine transformations. We will say that $\nabla_M$ is an affine connection.

Let $G$ be a Lie group, and $\hat G$ its universal cover, a left invariant affine structure defined on $G$, is an affine connection $\nabla_G$ invariant by the left translations of $G$. 

\medskip

The restriction of the connection to ${\cal G}$ defines a left-symmetric product,
 that is a product such that:

$$
ab-ba =[a,b].
$$

$$
a(bc)-(ab)c=b(ac)-(ba)c.
$$

Let $a$ be an element of ${\cal G}$, we denote by $L_a$, the linear map of ${\cal G}$ defined by $L_a(b)=ab$. We denote by $h(a)=(L_a,a)$. The following proposition has been shown by Koszul:

\bigskip

{\bf Proposition 3.4.}
{\it Let $G$ be a Lie group endowed with a left-symmetric algebra. The map $h(a)$ is an affine representation of the Lie algebra ${\cal G}$. Let $H:\hat G\rightarrow Aff({\cal G})$ be the representation whose differential at the origin is $h$. The orbit at the origin of $H$ is open. Conversely, let $\hat G$ be an $n$-dimensional simply connected Lie group, suppose that there exists an $n$-dimensional affine  representation $H$ of $G$ which has an open orbit, then ${\cal G}$ can be endowed with the structure of a left symmetric algebra which induces $H$.}

\bigskip

From this proposition, we deduce:

\bigskip

{\bf Proposition 3.5.}
{\it Suppose that $M$ is compact, and $Aff(M)_0$ acts transitively on $M$, then $M$ is a solvmanifold endowed with the structure of an affine manifold.}

\medskip

{\bf Proof.}
The corollary 3.1 shows that   $Aff({\bf R}^n,M)_0$ is solvable. Since $Aff({\bf R}^n,M)_0$ acts transitively and freely on ${\bf R}^n$, it is endowed with the structure of a left-invariant affine structure. We can identify $\pi_1(M)$ with a discrete subgroup of $Aff({\bf R}^n,M)_0$, since it action on $M$ is transitive and locally free. This implies that $M$ is an affine solvmanifold.

\bigskip

{\bf 4. $dim(Aff(M)_0)= dim(M)-1$.}

\bigskip

In this section, we are going to study compact ${\bf R}^n$-aspherical manifolds such that $dim(Aff(M)_0)=n-1$. In this case, the orbit of $Aff(M)_0$ are leaves of an $n-1$-dimensional foliations defined on $M$. Such foliations have been studied by many authors. Suppose that $Aff(M)_0$ is commutative, the results of  [2] imply:

\bigskip

{\bf Proposition 4.1.}
{\it Suppose that $M$ is an $n$-aspherical compact differentiable manifold such that $Aff(M)_0$ is commutative and has a compact orbit. Then $M$ is the total space of a fibre bundle over the circle $S^1$ with fibre the torus ${\bf T}^{n-1}$.}

\bigskip

In [4], Hector, Ghys and Moriyama have studied closed manifolds endowed with locally free actions of nilpotent Lie groups. The theorem $B$ of their paper  implies that if $Aff(M)_0$ is nilpotent, there exists a nilmanifold $N$ of dimension superior or equal to $1$ which is the fibre space of a locally trivial fibration whose total space is $M$. 

The corollary 3.4  implies: 

\bigskip

{\bf Proposition 4.2.}
{\it Suppose that $M$ is a compact differentiable $n$-aspherical manifold such that $Com({\bf R}^n, M)_0$ is an algebraic group and the dimension of its orbits is $n-1$. Then $M$ is the total space of a fibre space whose fibre is a nilmanifold of dimension superior or equal to $1$.}

\bigskip

The corollary 3.5 implies:

\bigskip

{\bf Corollary 4.1.}
{\it Let $M$ be an $n$-dimensional compact polynomial manifold such that $dim(Com({\bf R}^n,M)_0)=n-1$.  Then $M$ is the total space of a fibre space whose fibre is a nilmanifold of dimension superior or equal to $1$.}

\bigskip

In dimension $3$, we can precise these results: The solvable $2$-dimensional Lie groups are either isomorphic to ${\bf R}^2$, or to $GA$, the group of affine transformations of the real line. The action of $GA$ has been studied by Ghys
. We can deduce:

\bigskip

{\bf Theorem 4.2.}
{\it Let $M$ be a differentiable ${\bf R}^3$-aspherical manifold. Suppose that the dimension of $Aff(M)$ is $2$, and $Aff(M)_0$ is not commutative and preserves a volume form on $M$. Then $M$ is a solvmanifold, or the quotient of $\widetilde{SL(2,{\bf R})}$ by a lattice. }

\bigskip

\bigskip
\bigskip

\centerline{\bf Bibliography.}

\bigskip
\bigskip

1. Bredon, Glen E. Topology and geometry. Vol. 139. Springer Science  Business Media, 2013.

\smallskip

2. Chatelet, Gilles, Harold Rosenberg, and Daniel Weil. "A classification of the topological types of $\mathbf {R}^ 2$-actions on closed orientable 3-manifolds." Publications Mathématiques de l'IHÉS 43 (1974): 261-272.

\smallskip

3. Fried, David, William Goldman, and Morris W. Hirsch. "Affine manifolds with nilpotent holonomy." Commentarii Mathematici Helvetici 56, no. 1 (1981): 487-523.

\smallskip

4. Ghys, E. Hector, G. Moriyama, Y. On codimension $1$ nilmanifolds and a theorem of Malcev. Topology. 28 (1989): 197-210

\smallskip

5. Ghys, Etienne. "Action localement libres du groupe affine." Inventiones mathematicae 82, no. 3 (1985): 479-526.

\smallskip

6. Goldman, W. Hirsch, M. Affine manifolds and orbits of algebraic groups. Trans. Amer. Math. Soc. 295 (1986), 175-198.

\smallskip

7. Koszul, Jean-Louis. "Domaines bornés homogenes et orbites de groupes de transformations affines." Bulletin de la Société Mathématique de France 89 (1961): 515-533.

\smallskip

8. Poincar\'e, Henri. "Oeuvres de Henri Poincare Tome 1." (1928).

\smallskip

9. Rosenberg, H. Roussarie, R. Weil, D. A classification of closed $3$-manifolds of rank $2$. Ann. of Math., 91 (1970), 449-464.

\smallskip

 10. Tsemo, Aristide. "Dynamique des variétés affines." Journal of the London Mathematical Society 63, no. 2 (2001): 469-486.
\end{document}